\documentclass{amsart}
\usepackage[charter]{mathdesign}
\usepackage{palatino}
\usepackage{enumerate}
\theoremstyle{mythm}
\newtheorem{thm}{Theorem}[section]
\newtheorem{lem}[thm]{Lemma}
\newtheorem{prop}[thm]{Proposition}
\newtheorem{cor}[thm]{Corollary}
\theoremstyle{mydefn}

\theoremstyle{myrem}

\numberwithin{equation}{section}

\begin{document}
\title{More on cyclic amenability of the Lau product of Banach algebras defined by a Banach algebra morphism}
\author{M. Ramezanpour}
\address{School of Mathematics and Computer Science, 
Damghan University, P. O. Box 36716, Damghan 41167, Iran.}
\email{ramezanpour@du.ac.ir}
\date{}
\maketitle
\begin{abstract}
For two Banach algebras $A$ and $B$, the $T$-Lau product $A\times_T B$,
was recently introduced and studied for some bounded homomorphism $T:B\to A$ with $\|T\|\leq 1$. Here, 
we give general nessesary and sufficent conditions for $A\times_T B$ to be (approximately) cyclic amenable. 
In particular, we extend some recent results on (approximate) cyclic amenability of 
direct product $A\oplus B$ and $T$-Lau product $A\times_T B$ and answer a 
question on cyclic amenability of $A\times_T B$.
\end{abstract}
\section{Introduction and some Preliminaries}
The notion of weak amenability for commutative Banach algebras was introduced and studied for the
first time by Bade, Curtis and Dales \cite{B-C-D.awBLa}. 
Johnson \cite{John.DfL1iL1} extended this concept to the non commutative case and showed that group
algebras of all locally compact groups are weakly amenable. A Banach algebra $A$
is called weakly amenable if every continuous derivation $D:A\to A^*$ is inner.
It is often useful to restrict one's attention to derivations $D:A\to A^*$ satisfying the property
$D(a)(c)+D(c)(a)=0$ for all $a,c\in A$.
Such derivations are called cyclic. Clearly inner derivations are
cyclic. A Banach algebra is called cyclic amenable if every continuous cyclic derivations $D:A\to A^*$ 
is inner. This notion was presented by Gronbaek \cite{Gron.WcancBa}. He investigated
the hereditary properties of this concept, found some relations between cyclic
amenability of a Banach algebra and the trace extension property of its ideals. 

Ghahramani and Loy \cite{Ghah-Loy.Gna} introduced several approximate 
notions of amenability by requiring that all bounded derivations from a given Banach
algebra $A$ into certain Banach $A$-bimodules to be approximately inner. 
 In the same paper and the subsequent one \cite{Ghah-Loy-Zh.Gna2}, 
 the authors showed the distinction between each of these concepts and the corresponding
classical notions and investigated properties of algebras in each of these new classes. 
Motivated by this notions, Esslamzadeh and Shojaee \cite{Essl-Shoj.AwaBa} 
defined the concept of approximate cyclic amenability for
Banach algebras and investigated the hereditary properties for this new notion.
Shojaee and Bodaghi in \cite[Theorem 2.3]{Shoj-Bod.AgcaBa} showed that for Banach algebras $A$ 
and $B$, if direct product $A\oplus B$
with $\ell^1$-norm is approximately cyclic amenable, then so are $A$ and $B$. They also showed that the 
converse for the case where $A^2$ is dense in $A$.

On the other hand, for two Banach algebras $A$ and $B$ and a bounded homomorphism $T:B\to A$ with $\|T\|\leq 1$, 
the $T$-Lau product $A\times_T B$ is defined as the space $A\times B$ equipped with
the norm $\|(a,b)\|=\|a\|+\|b\|$ and the multiplication 
$$(a,b)(a',b')=(aa'+aT(b')+T(b)a', bb'),$$
for all $a, a'\in A$ and $b, b'\in B$. 
This product was first introduced and studied by
Bhatt and Dabhi in \cite{Bha-Dab.AraLpBadBam} for the case where $A$ is commutative.
Javanshiri and Nemati in \cite{Javn-Nemati.OcpBasp}
extended this product to the general Banach algebras and studied 
Arens regularity, amenability and n-weak amenability of $A\times_T B$ in the general case; 
see also \cite{Abt-Ghaf-Rej.BipBifLpBa} and \cite{Dab-Jab-Hagh.SnawaLpBadBam}.
When $T=0$, this multiplication is the usual coordinatewise product and so
$A\times_T B$ is in fact the direct product $A\oplus B$. Furthermore,
let $A$ be unital with the identity element $e$ and let $\theta: B\to \mathbb{C}$ be a non-zero multiplicative linear
functional. Define $T_\theta: B \to A$ as $T_\theta(b) = \theta(b)e$, for each $b\in B$. Then the above
product with respect to $T_\theta$ coincides with the product investigated in \cite{SM.laupro}. 

 Bhatt and Dabhi in \cite{Bha-Dab.AraLpBadBam} showed that 
 cyclic amenability of $A\times_T B$ is stable with respect to $T$, for the case where $A$ is commutative, but the proof  
 contains a gap. In \cite{Abt-Ghaf.AncaLpBadBam} Abtahi and Ghafarpanah  
fixed this gap and extended this result to an arbitrary Banach algebra $A$. Indeed 
they proved that if $A\times_T B$ is 
cyclic amenable then both $A$ and $B$ are cyclic amenable. 
They also proved the converse 
for the case where both $A$ and $B$ have faithful dual spaces. But they left it open for
all Banach algebras; \cite[Question 3.5]{Abt-Ghaf.AncaLpBadBam}. 

In the present paper, we give general nessesary and sufficent conditions 
for $A\times_T B$ to be (approximate)  cyclic amenable. 
In particular we extend the recent results on (approximate) cyclic amenability of the direct product 
$A\oplus B$, \cite{Shoj-Bod.AgcaBa}, and the $T$-Lau product $A\times_T B$ and answer 
Question 3.5 in \cite{Abt-Ghaf.AncaLpBadBam} on cyclic amenability of $A\times_T B$.

\section{Cyclic amenability}
Let $A$ be a Banach algebra, and  $X$ be a Banach $A$-bimodule. 
Then the dual space $X^*$ of  $X$ becomes a dual Banach $A$-bimodule with
the module actions defined by
$(fa)(x)=f(ax)$ and $(af)(x)=f(xa),$
for all $a\in A, x\in X$ and $f\in X^*$. 
A derivation from $A$ into $X$ is a linear mapping $D:A\to X$ satisfying
$$D(ac)= D(a)c + aD(c)\qquad(a, c\in A).$$
If $x\in X$ then $d_x:A\to X$ defined by $d_x(a)=ax-xa$ is a derivation which is called an inner derivation.

A derivation $D:A\to A^*$ is  said to be cyclic if $D(a)(c)+D(c)(a)=0$ for all $a,c\in A$.
If every continuous cyclic derivation $D:A \to A^*$ is inner then $A$ is  called cyclic amenable. 

As remarked in \cite{Bha-Dab.AraLpBadBam}, the dual space $(A\times_T B)^*$
can be identified with the Banach space $A^*\times B^*$ equipped with the norm  $\|(f,g)\|=\max\{\|f\|, \|g\|\}$ via
$$(f,g)(a,b)=f(a)+g(b),$$
where $a\in A, f\in A^*, b\in B $ and $g\in B^*$.
Moreover, a direct verification reveals that $(A\times_T B)$-module operations of $(A\times_T B)^*$
are as follows.
\begin{align*}
&(f,g)(a,b)=\left(fa+fT(b), T^{*}(fa)+gb\right),\\
&(a,b)(f,g)=\left(af+T(b)f, T^{*}(af)+bg\right)
\end{align*}
for $a\in A, b\in B, f\in A^{*}$ and $ g\in B^{*}$.

To clarify the relation between cyclic amenability of $A\times_T B$ and
that of $A$ and $B$, we need the next result which characterize the continuous cyclic derivations 
on $A\times_T B$.

\begin{lem}\label{Cder-T}
Suppose that $A$ and $B$ are Banach algebras and 
$T:B\to A$ is a bounded (by one) homomorphism. 
Then every  continuous cyclic derivation $D:A\times_T B\to (A\times_T B)^{*}$ enjoys the presentation
$$D(a,b)=(D_1(a)-S^{*}(b),D_2(b)+S(a))$$
for all $a\in A$ and $b\in B$, where 
\begin{enumerate}[\indent(a)]
\item $D_1:A\to A^{*}$ and $D_2:B\to B^{*}$ are continuous cyclic derivations.
\item $S:A\to B^{*}$ is a bounded linear operator such that 
$S(ac)=(T^{*}\circ D_1)(ac)$  and 
$S(a)b=(T^{*}\circ D_1)(a)b$   for all $a, c\in A$ and $b\in B$.
\end{enumerate}
Moreover, $D$ is inner if and only if  $S=T^{*}\circ D_1$ and both $D_1$ and $D_2$ are inner derivations. 
\end{lem}
Abtahi and Ghafarpanah in \cite{Abt-Ghaf.AncaLpBadBam}, proved that if $A\times_T B$ is 
cyclic amenable then both $A$ and $B$ are cyclic amenable. 
They also proved the converse 
for the case where both $A$ and $B$ have faithful dual spaces,  but left it  as an open question for
all Banach algebras; \cite[Question 3.5]{Abt-Ghaf.AncaLpBadBam}. An improvement of this result 
has been also obtained  in \cite{Nem-Jav.ScnATB}. Indeed, the converse has been proved for the case where
$A^2$ is dense in $A$; \cite[Theorem 2.6]{Nem-Jav.ScnATB}. We should remark that a Banach algebra $A$
has left (right) faithful dual space just when $A^2$ is dense in $A$.

Here we gives general necessary and sufficient conditions for $A\times_T B$ 
to be cyclic amenable. This result improves 
\cite[Theorem 2.6]{Nem-Jav.ScnATB} and answers also Question 3.5 in \cite{Abt-Ghaf.AncaLpBadBam}.
\begin{thm}\label{CA-T}
Suppose that $A$ and $B$ are Banach algebras and 
$T:B\to A$ is a bounded (by one) homomorphism. 
Then $A\times_T B$ is  cyclic amenable if and only if the following statements hold.
\begin{enumerate}
\item $A$ and $B$ are  cyclic amenable.
\item $A^2$ is dense in $A$ or $B^{2}$ is dense in $B$.
\end{enumerate}
\end{thm}

\section{Approximate cyclic amenability}
Recall from \cite{Ghah-Loy.Gna} that a derivation $D:A\to X$ is called approximately inner if there exists a net
$\{x_\alpha\}\subseteq X$ such that
$D=\lim_\alpha d_{x_\alpha}$ in the strong operator topology.
A Banach algebra $A$ is called  approximately cyclic amenable, if every
continuous cyclic derivation $D:A \to A^*$ is approximately inner.  The concepts
of approximate cyclic amenability was introduced and 
studied in \cite{Essl-Shoj.AwaBa}; see also \cite{Shoj-Bod.AgcaBa}.

The next result characterizes approximately inner derivations on $A\times_T B$.
\begin{prop}\label{AppCder-T}
Suppose that $A$ and $B$ are Banach algebras and 
$T:B\to A$ is a bounded (by one) homomorphism. 
Let $D:A\times_T B\to (A\times_T B)^{*}$ be a continuous cyclic derivation given by
$D(a,b)=(D_1(a)-S^{*}(b),D_2(b)+S(a))$
for all $a\in A$ and $b\in B$. Then $D$ is approximately inner
if and only if $S=T^*\circ D_1$ and both $D_1$ 
and $D_2$ are approximately inner.
\end{prop}
%
%
%

Applying Proposition \ref{AppCder-T} and using a similar argument used in the proof of Theorem \ref{CA-T} we can prove the next theorem.
This result improves \cite[Theorem 2.3]{Nem-Jav.ScnATB}.
\begin{thm}\label{AppCA-T}
Suppose that $A$ and $B$ are Banach algebras and 
$T:B\to A$ is a bounded (by one) homomorphism. 
Then $A\times_T B$ is  approximately cyclic amenable if and only if the following statements hold.
\begin{enumerate}
\item $A$ and $B$ are  approximately cyclic amenable.
\item $A^2$ is dense in $A$ or $B^{2}$ is dense in $B$.
\end{enumerate}
\end{thm}

Shojaee and Bodaghi in \cite[Theorem 2.3]{Shoj-Bod.AgcaBa} showed that for Banach algebras $A$ 
and $B$, if direct product $A\oplus B$
is approximately cyclic amenable, then so are $A$ and $B$. They also showed that the 
converse for the case where $A^2$ is dense in $A$. 
Applying Theorems \ref{CA-T} and \ref{AppCA-T} for $T=0$, we get the next result
which extends \cite[Theorem 2.3]{Shoj-Bod.AgcaBa}.
\begin{cor}
Let $A$ and $B$ be Banach algebras. Then the direct product $A\oplus B$ is 
(approximately) cyclic amenable if and only if the following statements hold.
\begin{enumerate}
\item $A$ and $B$ are  (approximately) cyclic amenable.
\item $A^2$ is dense in $A$ or $B^{2}$ is dense in $B$.
\end{enumerate}
\end{cor}

Let $A$ be unital and $\theta:B\to\mathbb{C}$ be a non-zero multiplicative linear functional. 
Define $T_\theta(b):=\theta(b)e$. Then $A\times_{T_\theta} B$ 
is the $\theta$-Lau product $A\times_\theta B$, \cite{SM.laupro}.  As a consequence of  
Theorems \ref{CA-T} and \ref{AppCA-T}, we have the next result. 
\begin{cor}
Let $A$ be unital and $\theta$ be a non-zero multiplicative linear functional on $B$. 
Then $A\times_\theta B$ is (approximately) cyclic amenable if and only if
both $A$ and $B$ are (approximately) cyclic amenable.
\end{cor}



\end{document}